\documentclass[11pt, leqno]{amsart}
\usepackage{amsthm, amsmath}
\usepackage{amssymb} 
\usepackage{amscd}
\usepackage[curve,matrix,arrow,frame,tips]{xy}
\usepackage{verbatim}
\usepackage{graphicx}
\newtheorem{thm}{Theorem}[section]
\newtheorem{prop}[thm]{Proposition}
\newtheorem{lemma}[thm]{Lemma}

\newtheorem{Definition}[thm]{Definition}
\newtheorem{Remarknumb}[thm]{Remark}

\newtheorem{Remark}[thm]{Remark}

\newcounter{ex}[section]

\def\thfill{\null\nobreak\hfill}

\def\endproof{\thfill\vbox{\hrule
  \hbox{\vrule\hbox to 5pt{\vbox to 5pt{\vfil}\hfil}\vrule}\hrule}}

\DeclareMathOperator{\Ker}{Ker}

\renewcommand{\to}{\longrightarrow}

\begin{document}

\title[ On Parahoric subgroups]{ On Parahoric subgroups}
\author[T. Haines]{T. Haines *}
\thanks{*Partially supported by  NSF Focused Research Grant 0554254.}
%\thanks{*Partially supported by NSF
%grant DMS99-70378 and by a Sloan Research Fellowship.}
\address{Department of
Mathematics\\
University of Maryland\\
College Park\\
MD 20742\\
USA}
\email{tjh@math.umd.edu}
\author[M. Rapoport]{M. Rapoport}
\address{Mathematisches Institut der Universit\"at Bonn,  
Beringstrasse 1\\ 53115 Bonn\\ Germany.}
\email{rapoport@math.uni-bonn.de}

\date{\today}

\maketitle

%\tableofcontents

 \bigskip

 {\bf Abstract}
  { \footnotesize We give the proofs of some simple facts on parahoric subgroups and on Iwahori Weyl groups used in \cite{H}, \cite{PR} and in \cite{R}.} 
\medskip


\bigskip
%\section*{Introduction}

\bigskip
Let $G$  be a connected reductive group over a
strictly henselian discretely valued field $L$. Kottwitz defines in \cite{Ko}
 a functorial surjective homomorphism
\begin{equation}
\kappa_G : G(L) \longrightarrow X^\ast (\hat Z (G)^I).
\end{equation}
Here $I = Gal (\bar L/L)$ denotes the absolute Galois group of
$L$.

Let ${\mathcal B}$ be the Bruhat-Tits building of the adjoint group
of $G$. Then $G(L)$ operates on ${\mathcal B}$.
\begin{Definition}{\rm
A  {\it parahoric subgroup} of  $G(L)$ is a subgroup of the form
\[
K_F = {\rm Fix} \, (F) \cap \Ker \, \kappa_G,
\]
for a facet $F$ of ${\mathcal B}$. An  {\it Iwahori subgroup} of
$G(L)$ is the parahoric subgroup associated to an alcove of $\mathcal B$.}
\end{Definition}

\begin{Remarknumb}
{\rm  If $F' = g F$, then
\[
K_{F'} = g K_F g^{-1}.
\]
In particular, since $G(L)$ acts transitively on the set of all
alcoves, all Iwahori subgroups are conjugate.}
\end{Remarknumb}

We shall see presently that this definition coincides with
the one of Bruhat and Tits \cite{BTII},  5.2.6. They associate to a facet $F$ in $\mathcal B$ a smooth group scheme $\mathcal G_F$  over ${\rm Spec}\ O_L$ with generic fiber $G$ and the open subgroup $\mathcal G_F^\circ$ of it with the same generic fiber and with connected special fiber, and define the parahoric subgroup attached to $F$ as $P_F^\circ=\mathcal G_F^\circ(O_L)$; from their definition it follows that $P_F^\circ \subseteq {\rm Fix}(F)$.  We denote by $G(L)_1$
the kernel of $\kappa_G$. 

For a facet $F$, let as above $P^\circ_F$ be the
associated parahoric subgroup in the BT-sense.  Our first 
goal is to prove the following proposition.

\begin{prop} \label{main_prop}  For any facet $F$ in $\mathcal B$, we have
 $P^\circ_F = K_F$.
\end{prop}

\begin{proof}
a) If $G=T$ is a torus, then $T(L)_1$ is the unique Iwahori
subgroup, cf.  Notes at the end of \cite{R}, ${\rm n^o} 1$. Hence the result follows
in this case.

b) If $G$ is semisimple and simply connected, then $G(L)_1 =
G(L)$. The assertion therefore follows from \cite{BTII}, 4.6.32, which also proves that $\mathcal G_F = \mathcal G_F^\circ$, where $\mathcal G_F$ is the group scheme defined in loc.~cit. 4.6.26.

%5.2.9.

c) Let $G$ be such that $G_{\rm der}$ is simply connected. Let $D
= G/G_{\rm der}$. We claim that there is a commutative diagram with exact rows
\[
\begin{array}{ccccccccc}
1& \to & \mathcal G_{{\rm der},F}(O_L) & \to & \mathcal G^\circ_F(O_L) & \to &  \mathcal D^\circ(O_L) & \to &1 \\[10pt]
&&\big\| & & \big\downarrow  &&\big\| &&\\[10pt]
1 & \to & K_{\rm der} & \to & K_F  &\to& D(L)_1 && .

\end{array}
\]
The bottom row involves parahoric subgroups associated in our sense to the facet $F$, and the top row involves those defined in \cite{BTII} and comes by restricting the exact sequence
\[
1 \longrightarrow G_{\rm der}(L) \longrightarrow G(L) \longrightarrow D(L) \longrightarrow 1
\]
to $ O_L$-points of the appropriate Bruhat-Tits group schemes (and more precisely in the case of $\mathcal D^\circ(O_L)$, the group $\mathcal D$ is the lft Neron model of $D_L$; cf. the notes at the end of \cite{R}).  The vertical equalities result from a) and b) above.  The vertical arrow is an inclusion (but we need to justify its existence, see below) making the entire diagram commutative.

Let us first {\em construct} the top row.  The key point is to show that the map $G(L) \to D(L)$ restricts canonically to a surjective map $\mathcal G^\circ_F(O_L) \to \mathcal D^\circ(O_L)$.  We shall derive this from the corresponding statement involving the lft Neron models of $D$ and of a maximal torus $T \subset G$.  Let $S$ denote a fixed maximal $L$-split torus in $G$, and define $T = {\rm Cent}_G(S)$, a maximal torus in 
$G$ since by Steinberg's theorem $G$ is quasisplit.  Also define $T_{\rm der} := G_{\rm der} \cap T = {\rm Cent}_{G_{\rm der}}(S_{\rm der})$, where $S_{\rm der} = (S \cap G_{\rm der})^\circ$, a maximal $L$-split torus in $G_{\rm der}$.

Consider the lft Neron models $\mathcal T$, $\mathcal T_{\rm der}$, and $\mathcal D$ associated to $T$, $T_{\rm der}$, and $D$, cf. \cite{BLR}.  The map $\mathcal T(O_L) \to \mathcal D(O_L)$ is surjective, and this implies that $\mathcal T^\circ(O_L) \to \mathcal D^\circ(O_L)$ is also surjective, cf. \cite{BLR}, 9.6, Lemma 2.  By \cite{BTII}, 4.6.3 and 4.6.7 we have decompositions $\mathcal G^\circ_F(O_L) = \mathcal T^\circ(O_L) \mathfrak U_F(O_L)$ and $\mathcal G_{{\rm der},F}(O_L) = \mathcal T_{\rm der}(O_L) \mathfrak U_F(O_L)$, where $\mathfrak U_F$ denotes the group generated by certain root-group $O_L$-schemes $\mathfrak U_{F,a}$, which depend on $F$; these all fix $F$.  These remarks show that $G(L) \to D(L)$ restricts to a map  $\mathcal G^\circ_F(O_L) \to \mathcal D^\circ(O_L)$, and also that the latter map is {\em surjective}.  

The kernel of the latter map is contained in the subgroup of $G_{\rm der}(L)$ which fixes $F$, hence is precisely $K_{\rm der} = \mathcal G_{{\rm der},F}(O_L)$.  This completes our discussion of the first row of the diagram above. 

%\footnote{Having established the existence of the exact sequence of $\mathcal O_L$-group schemes 
%$$
%1 \to \mathcal G_{{\rm der},F} \to \mathcal G^\circ_F \to \mathcal D^\circ \to 1
%$$
%we can give an alternate proof that $\mathcal G^\circ_F(O_L) \rightarrow \mathcal D^\circ(O_L)$ is surjective. Indeed, 
%$\mathcal G^\circ_F \rightarrow \mathcal D^\circ$ is a smooth epimorphism of smooth group schemes over $O_L$, which 
%is surjective on special fibers since its kernel is connected.}.

Since $G_{\rm der}$ is simply connected, $G(L)_1$ is the inverse image of $D(L)_1$ under the natural projection $G \to D$; hence $\mathcal G^\circ_F(O_L)$ belongs to $G(L)_1$ hence is contained in  $K_F$; this yields the inclusion which fits in the above diagram and makes it commutative.  A diagram chase then shows that the inclusion $\mathcal G^\circ_F(O_L) \to K_F$ is a bijection.

d) To treat the general case choose a $z$-extension $G' \to G$,
with kernel $Z$ where the derived group of $G'$ is simply
connected. Since $X_\ast (Z)_I=X^\ast (\hat Z^I)$ is torsion-free, the induced
sequence
\begin{equation}\label{exseq}
0 \longrightarrow X^\ast (\hat Z^I) \longrightarrow X^\ast (\hat
Z (G')^I) \longrightarrow X^\ast (\hat Z(G)^I) \longrightarrow 0
\end{equation}
is exact. We therefore obtain an exact sequence
\[
1 \longrightarrow Z(L)_1 \longrightarrow G'(L)_1 \longrightarrow
G(L)_1 \longrightarrow 1.
\]
It follows that we also have an exact sequence
\[
1 \longrightarrow Z(L)_1 \longrightarrow K'_F \longrightarrow K_F
\longrightarrow 1.
\]
As in c) one shows that $G'(L) \rightarrow G(L)$ maps $\mathcal G'^\circ_F(O_L) = P'^\circ_F$ {\em onto} $\mathcal G^\circ_F(O_L) = P^\circ_F$; in particular, since $P'^\circ_F = K'_F$ by c), we deduce $P^\circ_F \subset G(L)_1$ and thus $P^\circ_F \subseteq K_F$. 
The equality $P^\circ_F = K_F$
then follows from  $P'^\circ_F = K'_F$.
\end{proof}

\begin{Remarknumb} \label{main_prop_rem}
{\rm Proposition \ref{main_prop} makes sense and still holds true if $F$ is replaced with any bounded non-empty subset $\Omega \subset \mathcal B$ which is contained in an apartment.  Indeed, one can follow the same proof, making only the following adjustment in the proof of b) where $G= G_{\rm sc}$: although \cite{BTII} 4.6.32 is restricted to $\Omega$ contained in a facet, the equality  $\mathcal G^\circ_\Omega(O_L) = \mathcal G_\Omega(O_L)$ ($=$ fixer of $\Omega$) which we need requires only the connectedness of the group scheme $\mathcal T$ occurring in loc.cit., which holds here by loc.~cit.~4.6.1.} 
\end{Remarknumb}

\medskip

Let $F$ be a facet contained in the apartment associated to
the maximal split torus $S$. Let $T$ be the centralizer of $S$, a
torus since by Steinberg's theorem $G$ is quasisplit. Let $N$ be
the normalizer of $S$. Let $K_F$ be the parahoric subgroup
associated to $F$. Let
\[
\kappa_T : T(L) \longrightarrow X^\ast (\hat T^I) = X_\ast (T)_I
\]
be the Kottwitz homomorphism associated to $T$, and $T(L)_1$ its
kernel.
\begin{lemma} \label{T_cap_K=T_1}
\[
T(L) \cap K_F = T(L)_1.
\]
\end{lemma}

\begin{proof} By functoriality of the Kottwitz homomorphisms, we see $T(L)_1 \subset G(L)_1$.  
The elements of $\Ker \; \kappa_T$ act trivially on the
apartment associated to $S$, hence the inclusion "$\supset$" is
obvious. The converse follows from the fact that $T(L)_1$ equals
${\mathcal T}^\circ (O_L)$ where ${\mathcal T}^\circ$ is the identity
component of the lft Neron model of $T$ (cf. Notes at the end of \cite{R}, ${\rm n^o} 1$.) and the fact that ${\mathcal T}^\circ (O_L)$ is the
centralizer of $S$ in $K_F$, comp. \cite{BTII},  4.6.4 or \cite{L}, 6.3.
\end{proof}

\medskip

Let $K_0$ be the Iwahori subgroup associated to an alcove
contained in the apartment associated to $S$.
\begin{lemma}
\[
N(L) \cap K_0 =T(L)_1.
\]
\end{lemma}

\begin{proof} An  element of the $LHS$ acts trivially on the
apartment associated to $S$, hence is contained in $T(L) \cap K_0
= T(L)_1$ by the previous lemma.
\end{proof}
\begin{Definition}{\rm 
Let $S \subset T \subset N$ be as before (maximal split torus, contained in a maximal torus, contained in its normalizer). The {\it Iwahori-Weyl
group} associated to $S$ is
\[
\widetilde W = N (L) / T(L)_1.
\]
{\rm Let $ W_0= N (L) /T(L)$ (relative vector Weyl
group). We obtain an obvious exact sequence}
\begin{equation}\label{exactseq}
0 \longrightarrow X_\ast (T)_I \longrightarrow\widetilde W \longrightarrow
W_0 \longrightarrow 0 .
\end{equation}}
\end{Definition}
\begin{prop}\label{Bruhat}
Let $K_0$ be the Iwahori subgroup associated to an alcove
contained in the apartment associated to the maximal split torus
$S$. Then
\[
G(L) = K_0 . N(L). K_0
\]
and the map $K_0 n K_0 \mapsto n \in\widetilde W$ induces a bijection
\[
K_0 \backslash G(L) / K_0 \simeq\widetilde W.
\]
More generally, let $K$ resp. $K'$ be parahoric subgroups
associated to facets $F$ resp. $F'$ contained in the apartment
associated to $S$. Let
\[
\begin{array}{cccc}
&\widetilde W^K &=& (N(L) \cap K) / T(L)_1\ ,  \\
 &\widetilde W^{K'} &=& (N(L) \cap K') / T(L)_1\ .
\end{array}
\]
Then
\[
K \backslash G(L) / K' \simeq\widetilde W^K \backslash\widetilde W /\widetilde W^{K'}\ .
\]
\end{prop}

\begin{proof}
 To $F$ there are associated the subgroups
\[
U_F \subset P_F \subset G(L),
\]
cf. \cite{BTI}, 7.1 (where $P_F$ is denoted $\hat P_F$) or \cite{L},
8.8. By \cite{BTI}, 7.1.8 or \cite{L}, 8.10 we have
\[
P_F \subset U_F . N(L).
\]
Since $U_F$ is contained in the parahoric subgroup $P^\circ_F
\subset P_F$, the equality $G(L) = P_F . N(L). P_{F'}$ (\cite{BTI}, 7.4.15
resp. \cite{L}, 8.17) implies the first assertion
\[
G(L) = K_F. N(L).K_{F'}.
\]
To see the second assertion, we follow closely the proof of  \cite{BTI}, 7.3.4. Assume that $n,n' \in N(L)$ with
$n' \in P^\circ_F \,n\, P^\circ_{F'}$, i.e.
\[
n' . n^{-1} \in P^\circ_F \cdot P^\circ_{nF'}.
\]
We choose points $f \in F$ and $f' \in n \cdot F'$ with
\[
P^\circ_F = P^\circ_f, \quad P^\circ_{n.F'} = P^\circ_{f'}.
\]
We choose an order on the root system of $S$, with associated
vector chamber $D$ such that
\[
f \in f' +D.
\]
Then $U^-_f \cdot U_{f'}^- \subset U_{f'}^-$, in the notation of
\cite{L}, 8.8. Hence, by \cite{BTII}, 4.6.7 (comp. \cite{L}, 8.10),
\begin{eqnarray*}
P^\circ_f.P^\circ_{f'} &=& \left[ (N(L) \cap P^\circ_f) \cdot
U_{f}^+ . U_{f}^-
\right] \cdot \left[ U_{f'}^- U_{f'}^+ \cdot (N(L) \cap P^\circ_{f'})\right]\\
&=& (N(L) \cap P^\circ_f) \cdot U_{f}^+ . P_{f'}^\circ \\
&=& (N(L) \cap P^\circ_f) \cdot U_{f}^+ . U_{f'}^+. U_{f'}^-
(N(L) \cap P^\circ_{f'}).
\end{eqnarray*}
It follows that there exist $m_1 \in N(L) \cap P^\circ_f$ and
$m_2 \in N(L) \cap P^\circ_{f'}$ such that
\[
m_1 \cdot n' n^{-1} m_2 \in U^+ (L).U^- (L).
\]
{}From the usual Bruhat decomposition it follows that
\begin{eqnarray*}
m_1 \cdot n' n^{-1} m_2 &=& 1, \quad \quad{\rm i.e.,}\\
m_1 n' (n^{-1} m_2 n) &=& n \ .
\end{eqnarray*}
Since $m_1 \in N(L) \cap K_F$ and $n^{-1} m_2 n \in N(L) \cap
K_{F'}$, the last equality means that
\[
n \equiv n' \quad {\rm in} \quad\widetilde W^K \backslash\widetilde W /\widetilde W^{K'}.
\]
\end{proof}

\begin{Remarknumb} [Descent]  \label{descent_rem}
{\rm Let $\sigma$ denote an automorphism of $L$ having fixed field $L^\natural$ such that $L$ is the strict henselization of $L^\natural$.  Let us assume $G$ is defined over $L^\natural$; we may assume $S$, and hence $T$ and $N$, are likewise defined over $L^\natural$ (\cite{BTII}, 5.1.12).  Assume $F$ and $F'$ are $\sigma$-invariant facets in $\mathcal B$. Write $K(L^\natural) = K^\sigma$ and $K'(L^\natural) = K'^\sigma$.  We have a canonical bijection
\begin{equation}\label{invK}
 K(L^\natural) \backslash G(L^\natural) /K'(L^\natural) ~ \widetilde{\rightarrow} ~ [K\backslash G(L)/K']^\sigma.
\end{equation}
To prove that the map is surjective, we use the vanishing of $H^1(\langle \sigma \rangle, K)$ and $H^1(\langle \sigma \rangle, K')$.  To prove it is injective, we use the vanishing of $H^1(\langle \sigma \rangle, K \cap gK'g^{-1})$ for all $g \in G(L^\natural)$.  The vanishing statements hold because $K$, $K'$, and $K \cap gK'g^{-1}$ are the $O_L$-points of group schemes over $O_{L^\natural}$ with connected fibers, by virtue of Proposition \ref{main_prop} and Remark \ref{main_prop_rem}.

Next, note that,  using a similar cohomology vanishing argument, 
\begin{align*}
\widetilde{W}^\sigma &= N(L^\natural)/T(L^\natural) \cap T(L)_1 =: \widetilde{W}(L^\natural) \\
(\widetilde{W}^{K})^\sigma &= N(L^\natural) \cap K/T(L^\natural) \cap T(L)_1 =: \widetilde{W}^K(L^\natural). 
\end{align*}
Now suppose that $F$ and $F'$ are $\sigma$-invariant facets contained in the closure of a $\sigma$-invariant alcove in the apartment of $\mathcal B$ associated to $S$.  Then the canonical map
\begin{equation}\label{invW}
\widetilde{W}^{K}(L^\natural) \backslash \widetilde{W}(L^\natural) / \widetilde{W}^{K'}(L^\natural) ~ \rightarrow ~ [\widetilde{W}^K \backslash \widetilde W / \widetilde{W}^{K'}]^\sigma
\end{equation}
is bijective.  Indeed, note first that $\widetilde{W}^K$ and $\widetilde{W}^{K'}$ are parabolic subgroups of the quasi-Coxeter group $\widetilde{W}$ (see Lemma \ref{quasi_Coxeter_grp_lemma} below), and that any element $x \in \widetilde{W}$ has a {\em unique} expression in the form $wx_0w'$ where $w \in \widetilde{W}^K$, $w' \in \widetilde{W}^{K'}$, such that $x_0$ is the unique minimal-length element in $\widetilde{W}^K x_0 \widetilde{W}^{K'}$, and $wx_0$ is the unique minimal-length element in $wx_0\widetilde{W}^{K'}$.  Secondly, note that $\sigma$ preserves these parabolic subgroups as well as the quasi-Coxeter structure and therefore the Bruhat-order on $\widetilde{W}$.  These remarks imply the bijectivity just claimed.  Putting (\ref{invK}) and (\ref{invW}) together, we obtain a bijection
\[
K(L^\natural) \backslash G(L^\natural)/ K'(L^\natural) ~ \widetilde{\rightarrow} ~ 
\widetilde{W}^K(L^\natural) \backslash \widetilde{W}(L^\natural) / \widetilde{W}^{K'}(L^\natural).
\]
}
\end{Remarknumb} 

\begin{Remark} {\rm  We  now compare the Iwahori Weyl group with a variant of it  in  \cite{T}.

 In \cite{T}, p. 32, the following group is
introduced. Let $T(L)_b$ be the maximal bounded subgroup of
$T(L)$. The {\it affine Weyl group} in the sense of \cite{T}, p.32, associated to $S$ is  the
quotient $\widetilde W' = N(L)/T(L)_b$. We obtain a morphism of exact
sequences
\[
\begin{array}{ccccccccc}
0 & \longrightarrow & X_\ast (T)_I & \longrightarrow &\widetilde W & \longrightarrow &  W_0 & \longrightarrow  &0 \\[10pt]
&&\big\downarrow & & \big\downarrow & & \big\| && \\[10pt]
0 & \longrightarrow &\Lambda & \longrightarrow &
\widetilde W' & \longrightarrow & W_0 & \longrightarrow  &0\ .
\end{array}
\]

\medskip

\noindent Here $\Lambda=X_\ast (T)_I / torsion$. Of course,
\[
\Lambda = {\rm Hom} (X^\ast (T)^I, {\bf Z}) =
T(L)/T(L)_b\ ,
\]
(see \cite{Ko}, 7.2.)  
It follows that the natural homomorphism from $\widetilde W$ to $\widetilde W'$ is surjective, with finite kernel
isomorphic to $T(L)_b/T(L)_1$. 

The affine Weyl group in the sense of \cite{T} also appears in a kind of Bruhat decomposition, as follows. 
Let
\[
v_G : G(L) \longrightarrow X^\ast (\hat Z (G)^I)/torsion
\]
be derived from $\kappa_G$ in the obvious way. Let $C$ be an alcove in
the apartment corresponding to $S$. We consider the subgroup
\begin{equation}
\widetilde{K}_0 = {\rm Fix} (C) \cap \Ker \,v_G.
\end{equation}
Then $K_0$, the Iwahori subgroup corresponding to $C$,  is a normal subgroup of finite index in $\widetilde K_0$.  In fact $\Ker \, v_G$ is the group denoted $G^1$ in \cite{BTII} 4.2.16, and so by loc.~cit.~4.6.28, $\widetilde K_0$ is the group denoted there $\widehat{P}^1_C$, in other words the fixer of $C$ in $G^1$.  Using loc.~cit.~4.6.3, 4.6.7, we have $K_0 = T(L)_1 \,  \mathfrak U_{C}(O_L)$ and $\widetilde K_0 = T(L)_b \, \mathfrak U_C(O_L)$, where $\mathfrak U_C$ is the group generated by certain root-group $O_L$-schemes $\mathfrak U_{C,a}$ which fix $C$.  It follows that 
$\widetilde K_0 = T(L)_b K_0$, and  
\[
T(L)_b / T(L)_1 \simeq \widetilde K_0 / K_0.
\]

In \cite{T}, 3.3.1, the affine Weyl
group $\widetilde W'$ appears in the Bruhat decomposition with respect
to $\widetilde K_0$
\begin{equation} \label{Bruhat'}
\widetilde K_0 \backslash G(L) / \widetilde K_0 \simeq \widetilde W'.
\end{equation}
The group $T(L)_b/T(L_1)$ acts freely on both sides of the Bruhat decomposition
\[
K_0 \backslash G(L)/ K_0 = \widetilde{W}
\]
of Proposition \ref{Bruhat}, and we obtain (\ref{Bruhat'}) by taking the quotients.
}
\end{Remark}

\begin{Remark}
{\rm In \cite{BTII} the building of $G(L)$ (sometimes called the enlarged building $\mathcal B^1$) is also considered; it carries an action of $G(L)$.  There is an isomorphism $\mathcal B^1 = \mathcal B \times V_G$, where $V_G: = X_*(Z(G))_I \otimes \mathbb R$.  Given a bounded subset $\Omega$ in the apartment of $\mathcal B$ associated to $S$, there is a smooth group scheme $\widehat{\mathcal G}_\Omega$ whose generic fiber is $G$ 
and 
whose $O_L$-points $\widehat{\mathcal G}(O_L)$ is the subgroup of $G(L)$ fixing $\Omega \times V_G$, in other words the subgroup $\widehat{P}^1_\Omega$ in $G^1$ fixing 
$\Omega$.  We have $(\widehat{\mathcal G}_\Omega)^\circ =  \mathcal G^\circ_\Omega$.  Thus, the above discussion and Proposition \ref{main_prop} show that
$$
(\widehat{\mathcal G}_\Omega)^\circ(O_L) = \widehat{\mathcal G}_\Omega(O_L) \cap {\rm ker} \, \kappa_G.
$$ 
}
\end{Remark}

\medskip

\begin{prop}
Let $K$ be associated to $F$ as above. Let ${\mathcal G}^\circ_F$ be
the $O_L$-form with connected fibres of $G$ associated to $F$, and let
$\bar{\mathcal G}^\circ_F$ be its special fiber. Then $\widetilde W^K$ is
isomorphic to the Weyl group of $\bar{\mathcal G}^\circ_F$.
\end{prop}

\begin {proof} Let $\bar{\mathcal S}^\circ \subset \bar{\mathcal T}^\circ
\subset \bar{\mathcal G}^\circ_F$ be the tori associated to $S$ resp.
$T$, cf. \cite{BTII}, 4.6.4 or \cite{L}, 6.3. The natural projection
$\bar{\mathcal G}^\circ_F \to \bar{\mathcal G}^\circ_{F, {\rm red}} =
\bar{\mathcal G}^\circ_F/R_u (\bar{\mathcal G}^\circ_F)$ induces an
isomorphism
\[
\bar{\mathcal S}^\circ  \buildrel \sim \over \longrightarrow
\bar{\mathcal S}^\circ_{{\rm red}},
\]
where the index ,,red'' indicates the image group in $\bar{\mathcal
G}^\circ_{F, {\rm red}}$. Let $\bar W$ denote the Weyl group of
$\bar{\mathcal S}^\circ_{{\rm red}}$ and consider the natural homomorphism
\[
N(L) \cap K_F \longrightarrow \bar W.
\]
The surjectivity of this homomorphism follows from \cite{BTII},
4.6.13 or \cite{L}, 6.10. An element $n$ of the kernel centralizes
$\bar{\mathcal S}^\circ_{{\rm red}}$ and hence also $\bar{\mathcal S}^\circ$.
But then $n$ centralizes ${\mathcal S}$ because this can be checked
via the action of $n$ on $X_\ast ({\mathcal S})$. But the centralizer
of ${S}$ in $K_F$ is $T(L) \cap K_F$, cf. \cite{BTII}, 4.6.4 or \cite{L},
6.3. Hence $n \in T(L) \cap K_F$, which proves the claim.
\end{proof}

\medskip

We will give the Iwahori-Weyl group $\widetilde W$ the structure of a quasi-Coxeter group, that is, a semi-direct product of an abelian group with a Coxeter group.  Consider the real vector spaces $V = X_*(T)_I \otimes \mathbb R = X_*(S) \otimes \mathbb R$ and $V' := X_*(T_{\rm ad})_I \otimes \mathbb R$, where $T_{\rm ad}$ denotes the image of $T$ in the adjoint group $G_{\rm ad}$.  The relative roots $\Phi(G,S)$ for $S$ determine hyperplanes in $V$ 
(or $V'$), and the relative Weyl group $W_0$ can be identified with the group generated by the reflections through these hyperplanes.  The homomorphism
\[
T(L) \rightarrow X_*(T)_I \rightarrow V
\]
derived from $\kappa_T$ can be extended canonically to a group homomorphism
\[
\nu: N(L) \rightarrow V \rtimes \, W_0
\]
where $V \rtimes \, W_0$ is viewed as a group of affine-linear transformations of $V$; see \cite{T}, $\S 1$.  Using $\nu$, Tits defines in \cite{T} 1.4  the set of affine roots $\Phi_{\rm af}$, which can be viewed as a set of affine-linear functions on $V$ or on $V'$.  There is a unique reduced root system $\Sigma$ such that the affine roots $\Phi_{\rm af}$ consist of the functions on $V$ (or $V'$) of the form $y \mapsto \alpha(y) + k$ for $\alpha \in \, \Sigma$, $k \in \mathbb Z$.  The group generated by the reflections through the walls in $V$ (or in $V'$) through the hyperplanes coming from $\Phi_{\rm af}$ is the affine Weyl group $W_{\rm af}(\Sigma)$, which we also sometimes denote $W_a$.  The group $W_{\rm af}(\Sigma)$ can be given the structure of a Coxeter group, as follows.

The apartment $\mathcal A$ in $\mathcal B$ associated to $S$ is a torsor for $V'$; we identify it with $V'$ by fixing a special vertex $x$ in $\mathcal A$. 
%(we let $x_1$ denote a point in the apartment associated to $S$ in the extended building of $G(L)$, whose projection to $\mathcal B$ is $x$).  
  Assume $x$ belongs to the closure of the alcove $C \subset \mathcal A$.  Let ${\bf S}$ denote the reflections in $V'$ through the walls of $C$.  Then $(W_{\rm af}(\Sigma), {\bf S})$ is a Coxeter group.   Writing a superscript $x$ to designate subgroups of the affine linear transformations ${\rm Aff}(V')$ of $\mathcal A$ with $x$ as ``origin'', we have a semi-direct product
\[
^xW_{\rm af}(\Sigma) = Q^\vee \rtimes \, ^xW(\, \Sigma)\ ,
\]
where $Q^\vee = Q^\vee(\Sigma)$ is the coroot lattice for $\Sigma$ viewed as translations on $V'$, and where $^xW(\Sigma)$ is the finite Weyl group, isomorphic to $W_0$ and consisting of the elements in $^xW_{\rm af}(\Sigma)$ which fix $x$.  

The group $N(L)$ acts on the apartment $\mathcal A$ via the homomorphism $\nu$ when the latter is viewed as taking values in the group $^xW_{\rm af}(\Sigma)$ of affine-linear transformations on $\mathcal A$ endowed with $x$ as the ``origin''.  

Let $T_{\rm sc}$ resp.~$N_{\rm sc}$ be the inverse images of $T\cap G_{\rm der}$ resp.~$N\cap G_{\rm der}$ in $G_{\rm sc}$.  Let $S_{\rm sc}$ denote the split component of $T_{\rm sc}$. 
Then $S_{\rm sc}$ is a maximal split torus in $G_{\rm sc}$, and $T_{\rm sc}$ resp.~$N_{\rm sc}$  is its centralizer resp. normalizer.  
Let $\widetilde{W}_{\rm sc} =N_{\rm sc}(L)/T_{\rm sc}(L)_1$ be the Iwahori Weyl group of $G_{\rm sc}$. By \cite{BTII}, 5.2.10, this can be identified with $W_{\rm af}(\Sigma)$.  More precisely, note that the Kottwitz homomorphism for $G_{\rm sc}$ is trivial; also parahoric subgroups in $G_{\rm sc}(L)$ are simply the stabilizers or also the pointwise fixers of facets in $\mathcal B$. Let
$K_{{\rm sc},0}$ be the Iwahori subgroup of $G_{\rm sc}(L)$ associated to a chamber $C \subset \mathcal A$ and let ${\bf S}$ be the set of
reflections about the walls of $C$. Then (\cite{BTII}, 5.2.10)  
$(G_{\rm sc}(L), K_{{\rm sc},0}, N_{\rm sc}(L), {\bf S})$ is a double Tits system and $\nu: N_{\rm sc}(L) \rightarrow \, V' \rtimes \, ^xW(\Sigma)$ induces an isomorphism $\widetilde{W}_{\rm sc} \ = \, ^xW_{\rm af}(\Sigma)$.  In particular, $\widetilde{W}_{\rm sc}$ is a Coxeter group.   

\begin{prop}\label{semidirect}
Let $x \in\mathcal B$ be a special vertex in the apartment corresponding to the maximal split torus $S$, and let $K$ be the associated maximal parahoric subgroup of $G(L)$. The subgroup $\widetilde W^K$ projects isomorphically to the factor group $W_0$, and the exact sequence (\ref{exactseq}) presents $\widetilde W$ as a semi-direct product
$$
\widetilde W=X_\ast(T)_I\rtimes W_0\ .
$$
\end{prop}

\begin{proof}  In view of Lemma \ref{T_cap_K=T_1} we need only to see that $\widetilde W^K$ maps  surjectively  to $W_0$.  In the case of $G_{\rm sc}$, we have an isomorphism $\widetilde{W}_{\rm sc} = \, ^xW_{\rm af}(\Sigma)$ which sends $\widetilde{W}^{K_{\rm sc}}$ onto $\, ^xW(\Sigma)$, which is identified with $N_{\rm sc}(L)/T_{\rm sc}(L) = N(L)/T(L) = W_0$.  The desired surjectivity follows.
\end{proof}

\medskip

By a result of Borovoi \cite{Bo}, there is an exact sequence
\[
0 \rightarrow X_*(T_{\rm sc}) \rightarrow X_*(T) \rightarrow X^*(\hat{Z}(G)) \rightarrow 0.
\]
Since $X_*(T_{\rm sc})$ is an induced Galois module (\cite{BTII}, 4.4.16), $X_*(T_{\rm sc})_I$ is torsion-free and so the map $X_*(T_{\rm sc})_I \rightarrow X_*(T)_I$ is injective, and just as with (\ref{exseq}), 
\[
X_*(T)_I/X_*(T_{\rm sc})_I = X^*(\hat{Z}(G) ^I).
\]
Using Proposition \ref{semidirect} for both $\widetilde W$ and $\widetilde W_{\rm sc} \cong W_{\rm af}(\Sigma) =: W_a$, we deduce the following lemma.

\begin{lemma} \label{quasi_Coxeter_grp_lemma}
There is an exact sequence
\[
1 \to W_a \to\widetilde W \to X^\ast (\hat Z (G)^I ) \to 1.
\]
The subgroup $\Omega \subset \widetilde W$ consisting of the elements which preserve the alcove $C$ maps isomorphically onto $X^\ast (\hat Z (G)^I)$, defining a quasi-Coxeter structure on $\widetilde W$, 
\[ 
\widetilde W = W_a \rtimes \Omega. 
\]
\end{lemma}

\smallskip

Let $\nu'$ denote the composition of $\nu$ with the map $V \rightarrow V'$.  By \cite{T}, 1.7, $\Phi_{\rm af}$ is stable under $\nu'(N(L)) \subset {\rm Aff}(V')$, hence in particular under the group of translations $\nu'(T(L))$, which may be identified with $X_*(T_{\rm ad})_I$.  (Note that $X_*(T_{\rm ad})_I$ is torsion-free, by \cite{BTII}, 4.4.16.)

\begin{lemma}
There are natural inclusions
$$
Q^\vee \subset X_*(T_{\rm ad})_I \subset P^\vee\ ,
$$
where $P^\vee \subset V'$ is the set of coweights for $\Sigma$.  Moreover, $\nu'$ gives an identification
$$X_*(T_{sc})_I = Q^\vee.$$
\end{lemma}

\begin{proof}
The translations $\nu'(T(L))$ permute the special vertices in $\mathcal A$, and $P^\vee$ acts simply transitively on them; hence $\nu'(T(L)) \subset P^\vee$.  The equality $X_*(T_{\rm sc})_I = Q^\vee$ follows from the isomorphism $\widetilde W_{\rm sc} = \, ^xW_{\rm af}(\Sigma)$ deduced above from \cite{BTII}, 5.2.10, and this implies $Q^\vee \subset X_*(T_{\rm ad})_I$.
\end{proof}

\smallskip

\begin{Remark}{\rm The affine Weyl group associated to $S$ also appears in \cite{BTII} in a slightly different way. In \cite{BTII},
5.2.11 the following subgroup of $G (L)$ is introduced,
\begin{equation}
G(L)' = {\mbox{ {\it subgroup generated by all parahoric subgroups of }}}
\;G(L).
\end{equation}
In loc.cit. it is shown that this coincides with
\begin{equation}
G(L)' = T(L)_1 . \varphi (G_{\rm sc} (L)),
\end{equation}
where $\varphi : G_{\rm sc} \to G$ is the natural homomorphism
from the simply connected cover of $G_{{\rm der}}$.
\begin{lemma}
We have $G(L)' = G(L)_1$.
\end{lemma}

\begin{proof} The inclusion "$\subset$" is obvious. If $G = T$ is
a torus the equality is obvious, and so is the case when $G$ is
semisimple and simply connected. Next assume that $G_{\rm der}$ is
simply connected and let $D = G/G_{\rm der}$. We obtain a commutative
diagram
\[
\begin{array}{ccc}
G(L) & \longrightarrow & X^\ast (\hat Z (G)^I)\\[10pt]
\big\downarrow && \big\| \\[10pt]
D(L) &\longrightarrow & X^\ast (\hat D^I)\ .
\end{array}
\]
The exact sequence $1 \to T_{\rm der} \to T \to D \to 1$ induces a
morphism of exact sequences
\[
\begin{array}{ccccccccc}
1& \to & T_{\rm der} (L) & \to & T(L) & \to & D(L)& \to &1 \\[10pt]
&&\big\downarrow& & \big\downarrow \kappa_T &&\big\downarrow \kappa_D &&\\[10pt]
&&X_\ast (T_{\rm der})_I & \to & X_\ast (T)_I &\to& X_\ast (D)_I &\to&
0\ .
\end{array}
\]
All vertical maps are surjective. Now let $g \in G(L)_1$. Then by
a simple diagram chase there exists $t \in T(L)_1$ with the same
image in $D$ as $g$, and hence $g = tg' \in T(L)_1 \cdot \varphi
(G_{\rm der} (L)) = G(L)'$.

To treat the general case, choose a $z$-extension
\[
1 \to Z \to \widetilde G \buildrel \pi \over \longrightarrow G \to 1,
\]
with $\widetilde G_{\rm der}$ simply connected. We obtain a morphism of
exact sequences with surjective vertical maps
\[
\begin{array}{ccccccccc}
1& \to & Z (L) & \to & \widetilde G(L)& \to & G(L)& \to &1 \\[10pt]
&&\big\downarrow& & \big\downarrow  &&\big\downarrow  &&\\[10pt]
&&X^\ast (\hat Z^I) & \to & X^\ast (\hat Z (\widetilde G)^I) &\to&
X^\ast (\hat Z (G)^I) &\to& 0\ .
\end{array}
\]
The diagram shows that $G(L)_1 = \pi ( \widetilde G(L)_1)$. On the
other hand, the equality $G(L)' = \pi (\widetilde G(L)')$ is easy to
see. We conclude by the equality $\widetilde G (L)' = \widetilde G (L)_1$
which we already know.
\end{proof}

\bigskip

The other way of introducing the affine Weyl group is now as follows.
Let $N (L)' = N(L) \cap G(L)'$ and $W' = N(L)' / T(L)_1$. Let
$K_0$ be the Iwahori subgroup associated to an alcove $C$ in the
apartment corresponding to $S$ and let ${\bf S}$ be the set of
reflections about the walls of $C$. Then (\cite{BTII}, 5.2.12) $ (G
(L)', K_0, N(L)', {\bf S})$ is a double Tits system and $W'$ is
the affine Weyl group of the affine root system of $S$.  In particular, the natural homomorphism 
$W' \to W_{a}$ is an isomorphism. 
}
\end{Remark}

\bigskip

\noindent{\bf Acknowledgements:}  We thank E. Landvogt for his help. He pointed us  exactly to the places in  \cite{BTI, BTII},  resp. \cite{L}
needed in the proofs above.

\end{document}